\documentclass{amsart}
\usepackage{amssymb, amsmath, amscd}
\usepackage{latexsym}
\usepackage{graphicx}

\newcommand{\lo}{\longrightarrow}

\newcommand{\BC}{\Bbb C}
\newcommand{\eps}{\epsilon}
\newcommand{\del}{\delta}

\newcommand{\WAP}{\operatorname{WAP}}

\newcommand{\repn}{representation\,}
\newcommand{\repns}{representations\,}

\newtheorem{prop}{Proposition}[section]
\newtheorem{theo}{Theorem}[section]
\newtheorem{lemm}{Lemma}[section]
\newtheorem{cor}{Corollary}[section]

\newtheorem{rem}{Remark}[section]

\begin{document}

\title[restricted algebras]{restricted algebras on inverse semigroups I, representation theory}
\date{}
\author[M. Amini, A.R. Medghalchi]{Massoud Amini, Alireza Medghalchi}
\address{Department of Mathematics, Tarbiat Modarres University, P.O.Box
14115-175,\linebreak Tehran,Iran,mamini@modares.ac.ir \linebreak
\newline Department of Mathematics, Teacher Training University,
Tehran, Iran, medghalchi@saba.tmu.ac.ir} 
\thanks{This research was supported by Grant 510-2090
of Shahid Behshti University}

\begin{abstract}
The relation between representations and positive definite
functions is a key concept in harmonic analysis on topological
groups. Recently this relation has been studied on topological
groupoids. This is the first in a series of papers in which we
have investigated a similar relation on inverse semigroups. We use
a new concept of "restricted" representations and study the
restricted semigroup algebras and corresponding $C^*$-algebras.
\end{abstract}

\maketitle

\section{Introduction.} A continuous complex valued function $u:G\lo\BC$
on a locally compact Hausdorff group $G$ is called {\bf
positive definite} if for all positive integers $n$ and all
$c_1,\dots,c_n\in\BC$, and $x_1,\dots,x_n\in G$, we have
$$\sum_{i=1}^n \sum_{j=1}^n \bar c_ic_j u(x_i^{-1}x_j)\geq 0.$$
Positive definite functions on $G$ are automatically bounded [11].
When $G$ is Abelian, the positive definite functions are just the
Fourier-Stieltjes transforms of Radon measures on the Pontryagin
dual $\hat G$ of $G$ [10]. The significance of positive definite
functions is their relation with \repns of $G$. For each \repn
$\{\pi, \mathcal H_\pi\}$ of $G$ and each (unit) vector
$\xi\in\mathcal H_\pi$, the map
$$x\mapsto <\pi(x)\xi,\xi>$$
on $G$ is positive definite. Conversely each positive definite
function on $G$ is of this form [11].

Piere Eymard in [10] used positive definite functions to introduce
the Fourier and Fourier-Stieltjes algebras $A(G)$ and $B(G)$ of a
(not necessarily Abelian) locally compact group $G$ (see also
[12]). If $G$ is a locally compact Abelian group, $A(G)$ and
$B(G)$ are the ranges of the Fourier and Fourier-Stieltjes
transforms on $L^1(G)$ and $M(G)$, respectively. Dunkl and Ramirez
defined a subalgebra $R(S)$ of the algebra $\WAP(S)$ for a
Hausdorff locally compact commutative topological semigroup $S$ in
[9](see also [13],14]). For a locally compact Abelian group $G$,
$R(G)=M(\hat{G})^{\hat{}}$. In [15] A.T. Lau defined a subalgebra
$F(S)$ of $\WAP(S)$ for a topological $*$-semigroup $S$ and then
showed that if $S$ has an identity, then $F(S)=<P(S)>$, where
$P(S)$ is the set of all bounded positive definite functions on
$S$. When $S$ is commutative, we have $F(S)\subseteq R(S)$. The
authors studied the Fourier and Fourier Stieltjes algebras $A(S)$
and $B(S)$ of a unital topological foundation $*$-semigroup in
[5]. Also different versions of $A(G)$ and $B(G)$ are introduced
and studied for measured and topological groupoids $G$ in [18],
[19], [22], and [23].

This is the first in a series of papers whose ultimate goal is a
theory of Fourier algebras on inverse semigroups [1], [2]. There
are many technical difficulties when one tries to do things
similar to the group case. The first major obstacle is that the
regular representation of an inverse semigroup looses its
connection with positive definite functions. To avoid this
difficulty we decided to introduce a new concept of
representations. The new objects are called restricted
representations. The basic idea is that we require the
homomorphism property of representations to hold only for those
pairs of elements in the semigroup whose range and domain "match".
This is quite similar to what is done in the context of groupoids,
but the representation theory of groupoids is much more involved (
in [2] we investigate the relation between \repns of inverse
semigroups and their associated groupoids ).

This paper is organized as follows: In section 2 we introduce the
concept of restricted representations for inverse semigroups. In
section 3 we introduce the restricted semigroup algebra and study
its properties. In the last section the restricted versions of
full and reduced semigroup $C^*$-algebras are introduced and
studied.

\section{Preliminaries}
All over this paper, $S$
denotes a unital inverse semigroup with identity $1$. Let us remind that an
inverse semigroup $S$ is a discrete semigroup such that for each $s\in S$ there is a
unique element $s^*\in S$ such that
$$ss^*s=s,\quad s^*ss^*=s^*.$$
Then one can show that $s\mapsto s^*$ is an involution on $S$. The
set $E$ of idempotents of $S$  consists of elements the form
$ss^*,\, s\in S$. $E$ is a commutative sub semigroup of $S$. There
is a natural order $\leq$ on $E$ defined by $e\leq f$ if and only
if $ef=e$. We refer the reader to [21] for more details.

A $*$-\repn of $S$ is a pair $\{\pi,\mathcal H_\pi\}$ consisting
of a (possibly infinite dimensional) Hilbert space $\mathcal
H_\pi$ and a map $\pi:S\to\mathcal B(\mathcal H_\pi)$ satisfying
$$\pi(xy)=\pi(x)\pi(y),\, \pi(x^*)=\pi(x)^*\quad(x,y\in S),$$
that is a $*$-semigroup homomorphism from $S$ into the inverse
semigroup  of partial isometries on $\mathcal H_\pi$. We loosely
refer to $\pi$ as the \repn and it should be understood that there
is always a Hilbert space coming with $\pi$. Let
$\Sigma=\Sigma(S)$ be the family of all $*$-representations $\pi$
of $S$ with $$\|\pi\|:=\sup_{x\in S}\|\pi(x)\|\leq 1.$$ For $1\leq
p<\infty$, $\ell^p(S)$ is the Banach space of all complex valued
functions $f$ on $S$ satisfying
$$\|f\|_p :=\big(\sum_{x\in S} |f(x)|^p\big)^{\frac{1}{p}}<\infty .$$
For $p=\infty$, $\ell^\infty(S)$ consists of those $f$ with
$\|f\|_\infty  :=sup_{x\in S} |f(x)|<\infty$. Recall that $\ell
^1(S)$ is a Banach algebra with respect to the product
$$(f*g)(x)=\sum_{st=x} f(s)g(t) \quad
(f,g\in \ell ^1(S)),$$
and $\ell ^2(S)$ is a Hilbert space with inner product
$$<f,g>=\sum_{x\in S} f(x)\overline{g(x)}\quad (f,g\in\ell^2(S)).$$

Let also put
$$\check{f}(x)=f(x^*), \,  \tilde{f}(x)=\overline{f(x^*)},$$
for each $f\in \ell ^p (S)  \quad(1\leq p\leq \infty)$. We say
that $f$ is {\bf symmetric}, if $f=\tilde f$.

Next, following [16], we introduce the associated groupoid of an
inverse semigroup $S$. Given $x,y\in S$, the {\bf restricted
product} of $x,y$ is $xy$ if $x^*x=yy^*$, and undefined,
otherwise. The set $S$ with its restricted product forms a
groupoid [16,3.1.4] which is called the {\bf associated groupoid}
of $S$ and we denote it by $S_a$. If we adjoin a zero element $0$
to this groupoid, and put $0^*=0$, we get an inverse semigroup
$S_r$ with the multiplication rule
$$x\bullet y=\begin{cases}
xy & \text{if}\;\; x^*x=yy^* \\
0 & \text{otherwise}
\end{cases}\quad (x,y\in S\cup\{0\}),$$
which is called the {\bf  restricted semigroup} of $S$. A {\bf
restricted representation} $\{\pi,\mathcal H_\pi\}$ of $S$ is a
map $\pi:S\to\mathcal B(\mathcal H_\pi)$ such that
$\pi(x^*)=\pi(x)^*\quad (x\in S)$ and
$$\pi(x)\pi(y)=\begin{cases}
\pi(xy) & \text{if}\;\; x^*x=yy^* \\
0 & \text{otherwise}
\end{cases}\quad (x,y\in S).$$

Let $\Sigma_r=\Sigma_r(S)$ be the family of all restricted
representations $\pi$ of $S$ with $\|\pi\|=\sup_{x\in
S}\|\pi(x)\|\leq 1$.  It is not hard to guess that $\Sigma_r(S)$
should be related to $\Sigma(S_r)$. Let $\Sigma_0(S_r)$ be the set
of all $\pi\in \Sigma(S_r)$ with $\pi(0)=0$. Note that
$\Sigma_0(S_r)$ contains all cyclic representations of $S_r$. Now
it is clear that, via a canonical identification,
$\Sigma_r(S)=\Sigma_0(S_r)$.

One of the central concepts  in the analytic theory of inverse
semigroups (as in [4], [5],[6], [9], [17]) is the {\bf left
regular \repn} $\lambda:S\lo B(\ell^2(S))$ defined by
$$\lambda(x)\xi(y)=\begin{cases}
\xi(x^*y) & \text{if}\;\; xx^*\geq yy^* \\
0 & \text{otherwise}
\end{cases}\quad (\xi\in \ell^2(S), x,y\in S).$$
This lifts to a faithful \repn  of $\ell^1(S)$ [24]. We define the
{\bf restricted left regular representation} $\lambda_r:S\lo
B(\ell^2(S))$ by
$$\lambda_r(x)\xi(y)=\begin{cases}
\xi(x^*y) & \text{if}\;\; xx^*=yy^* \\
0 & \text{otherwise}
\end{cases}\quad (\xi\in \ell^2(S), x,y\in S).$$

Let us check that $\lambda_r\in \Sigma_r(S)$.  Given $x\in S$, and
$\xi,\eta\in\ell^2(S)$ we have
$\|\lambda_r(x)\xi\|_2^2=\sum\limits_{xx^*=yy^*}|\xi(x^*y)|^2\leq\sum\limits_{z\in
S}|\xi(z)|^2= \|\xi\|_2^2$, so $\|\lambda_r\|=\sup_{x\in
S}\|\lambda_r(x)\|\leq 1$. Also,
\begin{align*}
<\lambda_r(x^*)\xi,\eta> & = \sum_{y\in S}
\lambda_r(x^*)\xi(y)\overline{\eta(y)}\\
&=\sum_{x^*x=yy^*}\xi(xy)\overline{\eta(y)}\\
&=\sum_{xx^*=zz^*}\xi(z) \overline{\eta(x^*z)}=<\xi,\lambda_r(x)\eta>.
\end{align*}
The last two sums are equal, because, if $x^*x=yy^*$, then by taking $xy=z$, we have
$x^*z=x^*xy=y$, and $zz^*=xyy^*x^*=xx^*$. On the other hand, if $xx^*=zz^*$, then by taking
$x^*z=y$, we have, $yy^*=x^*zz^*x=x^*x$.

Finally for each $x,y,z\in S$, $\xi\in\ell^2(S)$,

$$\lambda_r(xy)\xi(z)=\begin{cases}
\xi(y^*x^*z) & \text{if}\;\; xyy^*x^*=zz^* \\
0 & \text{otherwise},
\end{cases}$$

where as

$$\lambda_r(x)\lambda_r(y)\xi(z)=\begin{cases}
\xi(y^*x^*z) & \text{if}\;\; xx^*=zz^*,\, yy^*=x^*zz^*x \\
0 & \text{otherwise},
\end{cases}$$

Now conditions $xx^*=zz^*$ and $yy^*=x^*zz^*x$  imply conditions
$x^*x=yy^*$ and $xyy^*x^*=zz^*$, so $\lambda_r(x)\lambda_r(y)$ is
equal to $\lambda_r(xy)$, if $xx^*=zz^*$, and is equal to $0$,
otherwise. Hence $\lambda_r\in \Sigma_r(S)$, as claimed.

Next, given $\xi,\eta$ in $\ell^2(S)$ or $\ell^1(S)$, we define
$$(\xi\bullet\eta)(x)=\sum\limits_{x^*x=yy^*}\xi(xy)\eta(y^*)\quad(x\in S).$$
Then it is easy to check that in both cases this is a convergent sum. Also we clearly have
$$<\lambda_r(x^*)\xi,\eta>=\xi\bullet\tilde{\eta}(x),$$
for each $x\in S$, and $\xi,\eta\in \ell^2(S)$.  To see the
relation between this new dot product with the original
convolution product on $\ell^1(S)$, it is useful to note that for
$\xi,\eta\in \ell^1(S)$, $\xi\bullet \eta$ could be equivalently
presented as
$$\xi\bullet \eta(x)=\sum_s\sum_{st=x, x^*x=t^*t} \xi(s)\eta(t)\quad(x\in S).$$
This in particular shows that if $\xi,\eta\geq 0$, then
$\|\xi\bullet\eta\|_1\leq\|\xi *\eta\|_1$, something which fails in general.

Similarly one can define the {\bf restricted right  regular
representation} $\rho_r$ of $S$ in $\ell^2(S)$ by
$$\rho_r(x)\xi(y)=\begin{cases}
\xi(yx) & \text{if}\; xx^*=y^*y \\
0 & \text{otherwise}
\end{cases}\quad (\xi\in \ell^2(S), x,y\in S),$$
and observe that $\rho_r\in\Sigma_r(S)$ and
$$<\rho_r(x)\xi,\eta>=\tilde{\eta}\bullet\xi(x)\quad (x\in S,
\xi,\eta\in\ell^2(S)).$$
Also we have
$<\tilde\rho_r(\varphi)\xi,\eta>=\varphi\bullet(\check{\xi}
\bullet\bar{\eta})(1)$, for each $\varphi\in\ell^1(S)$ and $\xi,\eta\in\ell^2(S)$,
where $1$ is the identity of $S$.
Indeed
\begin{align*}
<\tilde\rho_r(\varphi)\xi,\eta>&=\sum\limits_{y\in
S}(\tilde\rho_r(\varphi)\xi)(y)\bar{\eta}(y)=\sum\limits_y\sum\limits_z
\varphi(z)(\rho_r(z)\xi)(y)\bar{\eta}(y)\\
&=\sum\limits_z\varphi(z)<\rho_r(z)\xi,\eta>=\sum\limits_z\varphi(z)(\tilde{
\eta}\bullet\xi)(z)\\
&=\sum\limits_z\varphi(z)
 (\tilde{\eta}\bullet\xi)^{\check{}}(z^*)= \varphi\bullet
(\tilde{\eta}\bullet\xi)^{\check{}}(1)\\
&=\varphi\bullet
(\check{\xi}\bullet\bar{\eta})(1).
\end{align*}
These identities simplify further calculations. Therefore we are
led to consider the algebra $\ell^1(S)$ with respect to the dot
product $\bullet$ (instead of convolution product $*$). We  devote
the next section to the study this new algebra.

\section{Reduced semigroup algebra}

In this section we show that for an inverse semigroup $S$,
$(\ell^1(S),\bullet,\, \tilde{}\,)$ is a Banach $*$-algebra with
an approximate identity.

\begin{prop} The dot product is associative.
\end{prop}
{\bf Proof} Let $\xi,\eta,\theta\in \ell^1(S)$ and $x\in S$, then
$$(\xi\bullet\eta)\bullet\theta(x)=\sum_{x^*x=yy^*}(\xi\bullet
\eta)(xy)\theta(y^*)$$
If we put $z=xy$ then $x^*z=x^*xy=y$ and $zz^*=xyy^*x^*=xx^*$. Conversely
if $xx^*=zz^*$ then putting $y=x^*z$ we get $xy=xx^*z=z$ and $yy^*==x^*zz^*x=x^*x$. Hence the
above sum is equal to
$$\sum_{xx^*=zz^*}(\xi\bullet\eta)(z)\theta(z^*x)=\sum_{xx^*=zz^*}
\sum_{z^*z=vv^*}\xi(zu)\eta(v^*)\theta(z^*x).$$
On the other hand
\begin{align*}
\xi\bullet(\eta\bullet\theta)(x)
&=\sum_{x^*x=yy^*}\xi(xy)(\eta\bullet\theta)(y^*) \\
&=\sum_{x^*x=yy^*}\sum_{yy^*=uu^*} \xi(xy)\eta(y^*u)\theta(u^*)
\end{align*}
Given $x\in S$, assume that $z,v$ satisfy $xx^*=zz^*$ and $z^*z=vv^*$. Then put
$u=x^*z$ and $y=x^*zv$. We have $uu^*=x^*zz^*x=x^*x$ and
$$yy^*=x^*zvv^*z^*x=x^*zz^*zz^*x=x^*zz^*x=x^*x=uu^*.$$

Conversely, if $u,y$ satisfy $x^*x=uu^*$ and $uu^*=yy^*$, then put $z=xu$ and
$v=u^*y$. We have $zz^*=xuu^*x^*=xx^*$
and $vv^*=u^*yy^*u=u^*u$, and $z^*z=u^*x^*xu=u^*u$, so $vv^*=u^*u$.
Hence the two double sums which represent $(\xi\bullet\eta)\bullet\theta(x)$
and $\xi\bullet(\eta\bullet\theta)(x)$ are indeed the same.\qed

\begin{theo} Under the usual norm, $(\ell^1(S),\bullet, \tilde{}\,)$ is a Banach $*$-algebra.
\end{theo}
{\bf Proof} We need only to check that $(f\bullet g)\,\tilde{}
=\tilde g\bullet\tilde f$ and $\|f\bullet
g\|_1\leq\|f\|_1\|g\|_1$, for each $f,g\in \ell^1(S)$. Fix $f,g\in
\ell^1(S)$. For each $x\in S$,
\begin{align*}
(f\bullet g)\,\tilde{}\,(x)&=\sum_{xx^*=yy^*} \bar f(x^*y)\bar g(y^*) =\sum_{xx^*=yy^*}
\tilde f(y^*x)\tilde g(y)\\
&=\sum_{x^*x=zz^*} \tilde f(z^*)\tilde g(xz)=\tilde g\bullet\tilde f(x).
\end{align*}
Next  for $s,x\in S$ put $J_{s,x}=\{t\in S: st=x, x^*x=t^*t\}$ and
note that for each $s,x,y\in S$ with $x\neq y$ we have
$J_{s,x}\cap J_{s,y}=\emptyset$. This justifies the last
inequality of the following calculation
\begin{align*}
\|f\bullet g\|_1 & = \sum_x|\sum_s\sum_{t\in J_{s,x}}  f(s)g(t)| \leq\sum_x \sum_s
\sum_{t\in J_{s,x}}|f(s)||g(t)|\\
&=\sum_s |f(s)|\sum_x\sum_{t\in J_{s,x}}|g(t)|\leq\sum_s |f(s)| \sum_{t}|g(t)|=\|f\|_1\|g\|_1.
\end{align*}
\qed
\begin{rem}
One may define the dot product $\bullet$ with the sum running
through all elements $y$ satisfying $yy^*\leq x^*x$ and relate it
to the classical left regular representation $\lambda$ quite
similar to what we did with $\lambda_r$, but then the dot product
$\bullet$ won't be associative in general. On the other hand,
there is no connection between the usual convolution product $*$
on $\ell^1(S)$ and \repn $\lambda$. For these reasons, it seems
that the restricted version is inevitable if one insists to keep
the relation between left regular \repn and the multiplication on
the semigroup algebra. This relation is the key for our
computations all over the paper.
\end{rem}

We denote the above Banach algebra with $\ell^1_r(S)$ and call it the {\bf
restricted semigroup algebra}  of $S$.
Each $\pi\in\Sigma_r(S)$ lifts to an *-representation
$\tilde{\pi}$ of $\ell^1_r(S)$ via

$\quad\quad\quad\quad\quad\quad
\quad\quad\quad\tilde{\pi}(f)=\sum_{x\in S} f(x)\pi(x) \quad (f\in
\ell^1_r(S)).$

To see that $\tilde{\pi}$  is a $*$-representation of
$\ell_r^1(S)$, let $f,g\in\ell_r^1(S)$ and note that
\begin{align*}
\tilde\pi(f\bullet g) &=\sum_{x\in S}(f\bullet g)(x)\pi(x)\\
&=\sum_{x\in S}\sum_{x^*x=yy^*} f(xy)g(y^*)\pi(x),
\end{align*}
where as
\begin{align*}
\tilde\pi(f)\tilde\pi(g) &=\sum_{s\in S}\tilde\pi(f) g(s)\pi(s)\\
&=\sum_{s\in S}\sum_{t\in S} f(t)g(s)\pi(t)\pi(s)\\
&=\sum_{s\in S}\sum_{t^*t=ss^*} f(t)g(s)\pi(ts).
\end{align*}
Now the last sums of these presentations are converted to each other
via change of variables $x=ts, y=s^*$ and $t=xy, s=y^*$.

Note that $\ell_r^1(S)$ is not necessarily unital (even if $S$ is
unital). However  we show that it contains a (not necessarily
bounded) approximate identity consisting of finitely supported
functions.

\begin{lemm}
Given $y\in S$, $e\in E$,
$$\delta_y\bullet\delta_{e}=\begin{cases}
\delta_y & y^*y=e \\
0 & \text{otherwise}
\end{cases}
\quad \text{and}\quad\delta_e\bullet\delta_y=\begin{cases}
\delta_y & yy^*=e \\
0 & \text{otherwise.}
\end{cases}$$
\end{lemm}
{\bf Proof} We have
$\delta_y\bullet\delta_{e}(z)=\sum\limits_{z^*z=tt^*}\delta_y(zt)
\delta_{e}(t^*)=0$, unless $t=e$ and $zt=y$ and $z^*z=tt^*$. If these
equalities hold, then $z^*z=tt^*=t$, so $y=zt=zz^*z=z$.
Therefore $(\delta_y\bullet\delta_{e})(z)=0$, unless $z=y$. Now
$(\delta_y\bullet\delta_{e})(y)=\sum_{y^*y=tt^*}\delta_y(yt
)\delta_{e}(t^*)=0$, unless $t=e$, $yt=y$, and $y^*y=tt^*$.
These equalities imply that $y^*y=tt^*=e$.
Conversely if $y^*y=e$, then
$$\del_y\bullet
\del_{e}(y)=\sum_{y^*y=tt^*}\del_y(yt)\del_{e}(t^*)=\del
_y(ye)=\del_y(y)=1.$$

The other statement is proved similarly.\qed

Now for each finite subset $F=\{x_1,\dots,x_n\}$ of $S$ let us put
$$i(F)=\{e\in E: e=xx^* \,\text{or}\,\, x^*x, \,\text{for some}\,\, x\in F\},$$
which is clearly a finite subset of $E$. Define
$$e_F=\sum_{e\in i(F)} \delta_e.$$

\begin{lemm}
Let $F,G$ be finite subsets of $S$ and $e_F, e_G$ be as above, then

$(i)$ For each $F_0\subseteq F$  and each $s\in F_0$,
$e_F\bullet\delta_s=\delta_s\bullet e_F=\delta_s,$

$(ii)$ $e_F\bullet e_G=e_G\bullet  e_F=\sum_{e\in i(F)\cap i(G)}
\delta_e,$ in particular if $G\subseteq F$, then $e_F\bullet
e_G=e_G\bullet e_F=e_G$.

$(iii)$ For each $f=\sum_{i=1}^\infty f(s_i)\delta_{s_i}\in\ell_r^1(S)$,
$$f\bullet e_F=\sum_{s_i^*s_i\in i(F)} f(s_i)\delta_{s_i},$$
and
$$e_F\bullet f=\sum_{s_is_i^*\in i(F)} f(s_i)\delta_{s_i}.$$

$(iv)$ If $f\in\ell_r^1(S)$ and $supp(f)\subseteq F$, then $f\bullet e_F=e_F\bullet f=f$.

\end{lemm}
{\bf Proof} For each $s\in F_0$, there are unique $e,f\in i(F)$
such that $s^*s=e$ and $ss^*=f$. By Lemma 3.1, for each $g\in
i(F)$, $\delta_s\bullet\delta_g$ is equal to $\delta_s$, if $g=e$,
and is $0$, otherwise. A similar statement, with $f$ replaced by
$e$, holds for the multiplication by $\delta_s$ from right. This
shows $(i)$. Also the above lemma shows that for each $e,f\in E$,
$\delta_e\bullet\delta_f=\delta_f\bullet\delta_e$ is $\delta_e$,
if $e=f$, and is $0$, otherwise. Hence
$$e_F\bullet e_G=\sum_{e\in i(F),  f\in i(G)} \delta_e\bullet \delta_f
=\sum_{e\in i(F)\cap i(G)} \delta_e,$$
which proves the first statement of $(ii)$. In particular if
$G\subseteq F$, we get $e_F\bullet e_G=e_G\bullet e_F=e_G$. Next
let $f\in\ell_r^1(S)$ and choose any $s_i\in supp(f)$. If
$s_i^*s_i\notin i(F)$, then by above lemma, $\delta_{s_i}\bullet
\delta_e=0$, for each $e\in i(F)$, and so $\delta_{s_i}\bullet
e_F=0$. If $s_i^*s_i\in i(F)$, then there is a unique $e\in i(F)$
such that $s_i^*s_i=e$, so again by above lemma,
$\delta_{s_i}\bullet e_F=\delta_{s_i}\bullet
\delta_e=\delta_{s_i}$. This proves the first equality in $(iii)$.
The proof of the second equality is similar. Finally $(iv)$
follows from $(iii)$.\qed

\begin{prop}
The Banach algebra $\ell^1_r(S)$ has a
(not necessarily bounded) two sided approximate identity consisting of positive, symmetric
functions of finite support.
\end{prop}
{\bf Proof} For each finite subset $F$ of $S$ let $e_F$ be as
above. It is clear that $e_F$ is a positive, symmetric function of
finite support. Given $f\in l_r^1(S)$ we have $\sum\limits_{s\in
S}|f(s)|<\infty$, so there are at most countably many $s\in S$,
say $s_1,s_2,\dots$, for which $f(s)\neq 0$. Then given $\eps>0$,
there is $N\geq 1$ s.t. $\sum_{i=N+1}^\infty |f(s_i)|<\eps$. Put
$F_0=\{s_1,\dots, s_N\}$ and let $F\supseteq F_0$, then by the
first equality in part $(iii)$ of the above lemma,
\begin{align*}
\|f-f\bullet e_F\|_1 &=\sum_{s_i^*s_i\notin i(F)} |f(s_i)|
\leq\sum_{s_i^*s_i\notin i(F_0)} |f(s_i)|\\
&\leq\sum_{i\geq N+1} |f(s_i)|<\eps.
\end{align*}
Similarly we have $\|f-e_F\bullet f\|_1<\eps$.\qed

Note that with $f$, $F_0$ and $F$ as above, a similar argument
could show that $\|f-e_F\bullet f\bullet e_F\|_1<\eps$. The above
result looks more interesting when one recalls that if $S$ is not
unital, $\ell^1(S)$ may fail to have a bounded (or even unbounded)
approximate identity. This could be the case for $\ell^1(S_r)$,
but, in the light of the next proposition, what we have shown is
that $\ell^1(S_r)$ has an approximate identity, provided that we
identify two elements of $\ell^1(S_r)$ which agree at $0$.

\begin{prop} The  restriction map $\tau: \ell^1(S_r)\to \ell_r^1(S)$ is a
surjective contractive Banach algebra homomorphism whose kernel is $\mathbb C\delta_0$.
The quotient Banach algebra $\ell^1(S_r)/\mathbb C\delta_0$  is
isometrically isomorphic to $\ell_r^1(S)$.
\end{prop}
{\bf Proof} Recall that $S_r=S^0$ is a semigroup with respect  to
the restricted product. Let us denote the convolution product on
$\ell^1(S_r)$ by $\tilde *$, then
$$\delta_x\tilde *\delta_y =\delta_{x\bullet y}=\delta_x\bullet \delta_y\quad(x,y\in S_r),$$
where the second equality is trivial if $x=0$ or $y=0$, and could
be easily checked (similar to the proof of Lemma 3.1) when $x,y\in
S$. This shows that $\tau$ is a homomorphism. All the other
assertions are trivial, except that we have to check
$\|\tau(f)\|_1=\|f+\mathbb C\delta_0\|$, for each $f\in
\ell^1(S_r)$. But the right hand side is the infimum over all
$c\in \mathbb C$ of $\|f+c\delta_0\|_1$, which is clearly obtained
at $c=-f(0)$, since $ \|f+c\delta_0\|_1=\|\sum_s f(s)\delta_s
+(c+f(0))\delta_0\|_1=\sum_s |f(s)|+|c+f(0)|. $\qed

\section{restricted semigroup $C^*$-algebras}

We know that when $S$ is an inverse semigroup (not necessarily
unital), the left  regular \repn $\lambda$ lifts to a  faithful
representation $\tilde\lambda$ of $\ell^1(S)$ [24]. In particular
$\|f\|_{\lambda}=\|\tilde\lambda(f)\|$ defines a $C^*$-norm on
$\ell^1(S)$. The completion of $\ell^1(S)$ in this norm is  called
the {\bf reduced $C^*$-algebra} of $S$ and is denoted by
$C^*_{\lambda}(S)$. Now let $\Sigma=\Sigma(S)$ and for each $f\in
\ell^1(S)$  define
$$\|f\|_{\Sigma}=\sup\{ \|\tilde{\pi}(f)\|:\pi\in\Sigma=\Sigma(S)\},$$
where
$$\tilde\pi(f)=\sum_{x\in S} f(x)\pi(x)\quad(f\in\ell^1(S)).$$
Then for each irreducible representation
$\theta:C_{\lambda}^*(S)\lo
 B(\mathcal H_{\theta})$, $\tilde{\pi}=\theta\circ\tilde\lambda$ is
an irreducible representation of $\ell^1(S)$ and so
$\|f\|_{\lambda}\leq\|f\|_{\Sigma}$ $(f\in \ell^1(S_r))$ [9]. In
particular $\|\cdot\|_{\Sigma}$ is a $C^*$-norm. This is the
largest $C^*$-norm on $\ell^1(S)$ and the corresponding enveloping
$C^*$-algebra is denoted by $C^*(S)$ and is called the {\bf (full)
$C^*$-algebra} of $S$ [8]. (For full and reduced $C^*$-algebras on
topological groups see [7], [20], for topological semigroups see
[3],[15]).
 Also $\tilde\lambda$ extends uniquely to an
*-epimorphism $\tilde\lambda:C^*(S)\lo C^*_{\lambda}(S)$.

Now let us consider unital  inverse semigroup $S$ and its
associated inverse $0$-semigroup $S_r$. In this section we want to
explore the same ideas for $\ell_r^1(S)$. As all of the above
results are valid for the non unital inverse semigroups, we can
freely apply them to $S_r$.

Recall  that $\Sigma_r(S)=\Sigma_0(S_r)$. In particular, the
restricted left regular \repn $\lambda_r$ of $S$ corresponds to a
\repn of $S_r$ which vanishes at $0$. Indeed let $\Lambda$ dente
the left regular \repn of $S_r$. Consider the closed subspace
$\ell_0^2(S_r):=\{\xi\in\ell^2(S_r): \xi(0)=0\}$ of $\ell^2(S_r)$,
and let $P_0:\ell^2(S_r)\to \ell_0^2(S_r)$, \,$\xi\mapsto
\xi-\xi(0)\delta_0$ be the corresponding orthogonal projection,
then the restriction map is an isomorphism of Hilbert spaces from
$\ell_0^2(S_r)$ onto $\ell^2(S)$, and under this identification,
$\lambda_r(s)=\Lambda(s)P_0$. By the Wordingham's theorem [24] we
know that $\tilde\Lambda$ is a faithful \repn of $\ell^1(S)$.
Unfortunately the relation $\tilde\lambda_r=\tilde\Lambda(.) P_0$
does not necessarily imply that $\tilde\lambda_r$ is also
faithful, but a rather straightforward argument (even much easier
than that of [24]) shows that $\tilde\lambda_r$ is faithful.

\begin{lemm}
$\tilde\lambda_r$ is faithful.
\end{lemm}
{\bf Proof} Fix any  $f\in\ell_r^1(S)$ with
$\tilde\lambda_r(f)=0$. Let $u\in S$ and put $t=uu^*$ and
$\xi=\delta_{u^*}\in\ell^2(S)$, then
\begin{align*}
0=\tilde\lambda_r(f)\xi(t)&=\sum_{s\in S} f(s)(\lambda_r(s)\delta_{u^*})(uu^*)\\
&=\sum_{ss^*=uu^*} f(s)\delta_{u^*}(s^*uu^*)\\
&=\sum_{ss^*=uu^*, s^*uu^*=u^*} f(s)=f(u),
\end{align*}
so $f=0$.\qed

\begin{cor} The Banach algebra $\ell_r^1(S)$ is semi-simple.\qed
\end{cor}

The above lemma shows that
$\|f\|_{\lambda_r}:=\|\tilde\lambda_r(f)\|$ defines a $C^*$-norm
on $\ell_r^1(S)$. We call the completion $C_{\lambda_r}^*(S)$ of
$\ell_r^1(S)$ in this norm, the {\bf restricted reduced
$C^*$-algebra} of $S$. Next consider
$$\|f\|_{\Sigma_r}=\sup\{ \|\tilde{\pi}(f)\|:\pi\in\Sigma_r=\Sigma_r(S)\}
\quad (f\in \ell_r^1(S)),$$ then clearly
$\|f\|_{\lambda_r}\leq\|f\|_{\Sigma_r}$ and so
$\|\cdot\|_{\Sigma_r}$ is also a $C^*$-norm. Since $\tilde\pi$'s
with $\pi\in\Sigma_r(S)$ exhaust all the non degenerate
$*$-representations of $\ell_r^1(S)$, this is indeed the largest
$C^*$-norm on $\ell_r^1(S)$. We call the completion $C_{r}^*(S)$
of $\ell_r^1(S)$ in this norm, the {\bf restricted full
$C^*$-algebra} of $S$. As in the classical case, $\tilde\lambda_r$
extends uniquely to an *-epimorphism $\tilde\lambda_r:C_r^*(S)\lo
C^*_{\lambda_r}(S)$.

Next we find the relation between the restricted full and reduced
$C^*$-algebras of $S$ with the full and reduced $C^*$-algebras of
$S_r$. The following technical lemma is probably true in a more
general form.

\begin{lemm}
Let $(A, \|.\|_1)$ be a Banach algebra and $\|.\|$ be a $C^*$-norm
on $A$ satisfying $\|.\|\leq \|.\|_1$. Let $C^*(A)$ be the
completion of $A$ in $\|.\|$. Let $J$ be a two sided ideal of $A$
which is closed in $(C^*(A),\|.\|)$. Then the quotient norm on
$C^*(A)/J$ induces a $C^*$-norm on $A/J$ and the $C^*$-completion
$C^*(A/J)$ of $A/J$ in this norm is isometrically isomorphic to
$C^*(A)/J$.
\end {lemm}
{\bf Proof} $J$ is clearly a closed ideal of $(A,\|.\|_1)$. Hence
$A/J$ is a Banach algebra under the quotient norm induced by
$\|.\|_1$ which is clearly dense in $C^*(A)/J$ in its
$C^*$-quotient norm. \qed

We apply the above lemma to $A=\ell^1(S_r)$, $J=\mathbb C\delta_0$, and $C^*$-norms
$\|.\|_{\Lambda}$ and $\|.\|_{\Sigma(S_r)}$. Note that $\mathbb C\delta_0$ is closed
in both of the above $C^*$-norms. Indeed, given $c\in \mathbb C$, we have
\begin{align*}
\|c\delta_0\|_{\sigma(S_r)} &=sup_{\pi\in\Sigma(S_r)}\|\tilde\pi(c\delta_0)\|
=sup_{\pi\in\Sigma(S_r)}\|c\pi(0)\|=|c|,\\
\end{align*}
and
\begin{align*}
\|c\delta_0\|_{\Lambda}&=\|\tilde\Lambda(c\delta_0)\|=\|c\Lambda(0)\|=|c|,
\end{align*}
where  the last equality follows from the fact that
$\|\Lambda(0)\|=1$ (inequality in one direction is already known,
for the other direction use $\Lambda(0)\delta_0=\delta(0)$).
Therefore any net $\{c_\alpha\delta_0\}$ which is convergent in
any of the above norms would result in a Cauchy net $\{c_\alpha\}$
in $\mathbb C$. If $c_\alpha\to c$ in $\mathbb C$, then the given
net would converge to $c\delta_0$ in both norms.

Therefore,  in the light of  Proposition 3.3, The following result
follows immediately from the above lemma.

\begin{prop} We have the isometric isomorphisms of $C^*$-algebras
$C_{\lambda_r}^*(S)\simeq  C_{\Lambda}^*(S_r)/\mathbb C\delta_0$
and $C_{r}^*(S)\simeq C^*(S_r)/\mathbb C\delta_0$.\qed
\end{prop}

{\bf acknowledgement.} We would like to thank  Dr. David Cowan who
pointed out an error in the first version of this paper. The first
author would like to thank hospitality of Professor Mahmood
Khoshkam during his stay in University of Saskatchewan, where the
main part of the revision was done.


\begin{thebibliography}{99}

\bibitem{} M. Amini, A. Medghalchi, restricted algebras on inverse semigroups II,
positive definite functions, preprint,
Shahid Beheshti University, 2000.

\bibitem{} M. Amini, A. Medghalchi, restricted algebras on inverse semigroups III,
Fourier algebra, preprint,
Shahid Beheshti University, 2000.


\bibitem{} M. Amini, A. Medghalchi, Fourier algebras on topological foundation
$*$-semigroups, preprint,
Shahid Beheshti University, 2000.

\bibitem{a} B. A. Barnes, Representations of the $l^1$-algebra of an
inverse semigroup, Trans. Amer. Math. Soc. 218 (1976) 361-396.

\bibitem{} C. Berg, J. P.R. Christensen, P. Ressel, Harmonic Analysis on
semigroups, Springer-Verlag, Berlin, 1984.

\bibitem{} J. B. Conway, J. Duncan, A. L. T. Paterson,
 Monogenic inverse semigroups and their $C^*$-algebras, Proc. Roy.
Soc. Edinburgh 98A (1984) 13-24.

\bibitem{} J. Dixmier, $C\sp*$-algebras, North-Holland Mathematical Library, Vol. 15,
North-Holland, Amsterdam, 1977.

\bibitem{a} J. Duncan, A. L. T. Paterson,  $C^*$-algebras of inverse
semigroups, Proc. Edinburgh Math. Soc. 28 (1985) 41-58.

\bibitem{a} C.F. Dunkl, D.E. Rumirez, $L^\infty$-representations of
commutative semitopological semigroups, Semigroup Forum 7 (1974) 180-199.

\bibitem{a} P. Eymard,  L'algebra de Fourier d'un groupe localement
compact, Bull. Soc. Math. France, 92 (1964) 181-236.

\bibitem{a} R. Godement,  Les fonctions de type positive et la theorie des
groupes, Trans. Amer. Math. Soc. 63 (1948) 1-84.

\bibitem{a} E. Hewitt, K.A. Ross, Abstract Harmonic Analysis I, second ed.,
Grundlehren der Mathematischen Wissenschaften 115, Springer-Verlag, Berlin, 1963.

\bibitem{a} M. Lashkarizadeh Bami,  Bochner's theorem and the Hausdorff
moment theorem on foundation topological semigroups, Can. J. Math. 37 (1985) 785-809.

\bibitem{a} M. Lashkarizadeh Bami,  Representations of foundation
semigroups and their algebras, Can. J. Math. 37 (1985) 29-47.

\bibitem{a} A. T. M. Lau,  The Fourier Stieltjes algebra of a topological
semigroup with involution, Pac. J. Math. 77 (1978) 165-181.

\bibitem{} Mark V. Lawson, Inverse semigroups, the theory of partial symmetries,
World Scientific, Singapore, 1998.

\bibitem{} R.J. Lindahl, P.H. Maserick, Positive-definite functions on involution
semigroups, Duke Math. J. 38 (1971) 771-782.

\bibitem{} K. Oty, Fourier-Stieltjes algebras of r-discrete groupoids,
J. Operator Theory, 41 (1999) 175-197.

\bibitem{} A.T. Paterson, The Fourier algebra for locally compact groupoids, preprint, 2002.

\bibitem{} G. K. Pedersen,  $C^*$-algebras and their automorphism groups,
Academic Press, New York, 1979.

\bibitem{} M. Petrich, Inverse semigroups, John Wiley, New York, 1984.

\bibitem{}
A. Ramsay, M.E. Walter, Fourier-Stieltjes algebras of locally compact groupoids,
J. Functional Analysis, 148 (1997) 314-365.

\bibitem{}
J. N. Renault, The Fourier algebra of a measured groupoid and its multipliers,
J. Functional Analysis, 145 (1997) 455-490.

\bibitem{a} J.R. Wordingham, The left regular \repn of an inverse semigroup,
Proc. amer. Math. Soc. 86   (1982) 55-58.

\end{thebibliography}
\end{document}